\documentstyle[11pt]{article}

\newtheorem{lemma}{{\bf Lemma}}[section]
\newtheorem{theo}{{\bf Theorem}}

\newtheorem{cor}[theo]{{\bf Corollary}}
\newtheorem{defn}{{\bf Definition}}[section]
\newtheorem{remark}{{\bf Remark}}[section]

\font\bbb=msbm10 scaled\magstep1

\newcommand{\RR}{\mbox{\bbb R}}
\newcommand{\ZZ}{\mbox{\bbb Z}}

\newcommand{\TPSS}{S^{\hspace{.1mm}1}\! \times \hspace{-3.3mm}_{-} \,
S^{\hspace{.2mm}2}}

\newcommand{\TPSSDS}{S^{\hspace{.1mm}1}\! \times \hspace{-3.3mm}_{-} \, S^{\hspace{.2mm}d-1}}

\newcommand{\La}{\Lambda}

\begin{document}

\title{\bf Minimal triangulations of
sphere bundles over the circle}
\author{{\bf Bhaskar Bagchi}$^{\rm a}$ and {\bf Basudeb Datta}$^{\rm b}$
}

\date{}

\maketitle

\vspace{-5mm}

\noindent {\small $^{\rm a}$Theoretical Statistics and
Mathematics Unit, Indian Statistical Institute,  Bangalore
560\,059, India

\smallskip

\noindent $^{\rm b}$Department of Mathematics, Indian Institute
of Science, Bangalore 560\,012,  India}

\footnotetext[1]{{\em E-mail addresses:} bbagchi@isibang.ac.in
(B. Bagchi), dattab@math.iisc.ernet.in (B. Datta). }

\begin{center}

%\date{December 8, 2006}
%\date{December 11, 2006}
%\date{Revised: August 20, 2007}
\date{To appear in `Journal of Combinatorial Theory, Ser. A'}

\end{center}

\hrule

\medskip

 \noindent {\bf Abstract}

\medskip

{\small For integers $d \geq 2$ and $\varepsilon = 0$ or $1$, let
$S^{\hspace{.1mm}1, d - 1}(\varepsilon)$ denote the sphere
product $S^{\hspace{.1mm}1} \times S^{\hspace{.1mm}d - 1}$ if
$\varepsilon = 0$ and the twisted sphere product $\TPSSDS$ if
$\varepsilon = 1$. The main results of this paper are\,: $(a)$ if
$d \equiv \varepsilon$ (mod 2) then $S^{\hspace{.1mm}1, d -
1}(\varepsilon)$ has a unique minimal triangulation using $2d+3$
vertices, and $(b)$ if $d \equiv 1 - \varepsilon$ (mod 2) then
$S^{\hspace{.1mm}1, d - 1}(\varepsilon)$ has minimal
triangulations (not unique) using $2d+4$ vertices. In this
context, a minimal triangulation of a manifold is a triangulation
using the least possible number of vertices. The second result
confirms a recent conjecture of Lutz. The first result provides
the first known infinite family of closed manifolds (other than
spheres) for which the minimal triangulation is unique. Actually,
we show that while $S^{\hspace{.1mm}1, d - 1}(\varepsilon)$ has
at most one $(2d + 3)$-vertex triangulation (one if $d \equiv
\varepsilon$ (mod 2), zero otherwise), in sharp contrast, the
number of non-isomorphic $(2d + 4)$-vertex triangulations of
these $d$-manifolds grows exponentially with $d$ for either
choice of $\varepsilon$. The result in $(a)$, as well as the
minimality part in $(b)$, is a consequence of the following
result\,: $(c)$ for $d\geq 3$, there is a unique $(2d +
3)$-vertex simplicial complex which triangulates a non-simply
connected closed manifold of dimension $d$. This amazing
simplicial complex was first constructed by K\"{u}hnel in 1986.
Generalizing a 1987 result of Brehm and K\"{u}hnel, we prove that
$(d)$ any triangulation of a non-simply connected closed
$d$-manifold requires at least $2d + 3$ vertices. The result
$(c)$ completely describes the case of equality in $(d)$. The
proofs rest on the Lower Bound Theorem for normal pseudomanifolds
and on a combinatorial version of Alexander duality.}

\medskip

{\small


\smallskip


} \medskip

\hrule

\section{Preliminaries}

With a single exception in Section 3, all simplicial complexes
considered here are finite. For a simplicial complex $X$, $V(X)$
will denote the set of all the vertices of $X$ and $|X|$ will
denote the geometric carrier of $X$. One says that $X$ is a {\em
triangulation} of the topological space $|X|$. If $|X|$ is a
manifold then we say that $X$ is a {\em triangulated manifold}.
The unique $(d + 2)$-vertex triangulation of the $d$-sphere
$S^{\hspace{.1mm}d}$ is denoted by $S^{\hspace{.1mm}d}_{d + 2}$
and is called the {\em standard $d$-sphere}. The unique
$(d+1)$-vertex triangulation of the $d$-ball is denoted by
$B^{\hspace{.1mm}d}_{d + 1}$ and is called the {\em standard
$d$-ball}. For $n\geq 3$, the unique $n$-vertex triangulation of
the circle $S^1$ is denoted by $S^1_n$ and is called the {\em
$n$-cycle}.

For $i = 1, 2$, the $i$-faces of a simplicial complex $K$ are
also called the {\em edges} and {\em triangles} of $K$,
respectively. For a simplicial complex $K$, the graph whose
vertices and edges are the vertices and edges of $K$ is called
the {\em edge graph} (or {\em $1$-skeleton}) of $K$. Recall that
a {\em graph} is nothing but a simplicial complex of dimension at
most $1$. A set of vertices in a graph is called a {\em clique}
if these vertices are mutually adjacent (i.e., any two of them
form an edge). Note that any simplex in a simplicial complex is a
clique in its edge graph.

For a simplex $\sigma$ in a simplicial complex $K$, the number of
vertices in ${\rm lk}_K(\sigma)$ is called the {\em degree} of
$\sigma$ in $K$ and is denoted by $\deg_K(\sigma)$ (or by
$\deg(\sigma)$). So, the degree of a vertex $v$ in $K$ is the same
as the degree of $v$ in the edge graph of $K$.

Recall that for any face $\alpha$ of a complex $X$, its {\em
link} ${\rm lk}_X(\alpha)$ is the simplicial complex whose faces
are the faces $\beta$ of $X$ such that $\alpha\cap \beta =
\emptyset$ and $\alpha\cup \beta\in X$. Likewise, the {\em star}
${\rm st}_X(\alpha)$ of the face $\alpha$ has all the maximal
faces $\gamma \supseteq \alpha$ of $X$ as its maximal faces.

A simplicial complex $X$ is called a {\em combinatorial
$d$-sphere} (respectively, {\em combinatorial $d$-ball\,}) if
$|X|$ (with the induced pl structure from $X$) is pl homeomorphic
to $|S^{\hspace{.2mm}d}_{d + 2}|$ (respectively,
$|B^{\hspace{.1mm}d}_{d + 1}|$). A simplicial complex $X$ is said
to be a {\em combinatorial $d$-manifold} if $|X|$ (with the
induced pl structure) is a pl $d$-manifold. Equivalently, $X$ is
a combinatorial $d$-manifold if all its vertex links are
combinatorial spheres or combinatorial balls. In this case, we
also say that $X$ is a {\em combinatorial triangulation} of
$|X|$. A simplicial complex $X$ is a combinatorial manifold
without boundary if all its vertex links are combinatorial
spheres. A combinatorial manifold will usually mean one without
boundary.

A simplicial complex $K$ is called {\em pure} if all the maximal
faces ({\em facets}) of $K$ have the same dimension.
%A maximal face in a pure complex is also called a {\em facet}.
For $d \geq 1$, a $d$-dimensional pure simplicial complex is said
to be a {\em weak pseudomanifold} if each $(d - 1)$-simplex is in
exactly two facets. Clearly, any $d$-dimensional weak
pseudomanifold has at least $d + 2$ vertices, with equality only
for $S^{\hspace{.2mm}d}_{d+2}$.

For a pure $d$-dimensional simplicial complex $K$, let $\La (K)$
be the graph whose vertices are the facets of $K$, two such
vertices being adjacent in $\La (K)$ if the corresponding facets
intersect in a $(d - 1)$-face. If $\Lambda(K)$ is connected then
$K$ is called {\em strongly connected}. A strongly connected weak
pseudomanifold is called a {\em pseudomanifold}. Thus, for a
$d$-pseudomanifold $K$, $\Lambda(K)$ is a connected $(d +
1)$-regular graph. This implies that $K$ has no proper subcomplex
which is also a $d$-pseudomanifold. (Or else, the facets of such
a subcomplex would provide a disconnection of $\Lambda(X)$.) By
convention, $S^{\hspace{.2mm}0}_2$ is the only $0$-pseudomanifold.

%\begin{defn}$\!\!\!${\bf .} \label{NPMFD} {\rm
A connected $d$-dimensional weak pseudomanifold is said to
be a {\em normal pseudomanifold} if the links of all the
simplices of dimension up to $d - 2$ are connected. Clearly, the
normal $2$-pseudomanifolds are just the connected combinatorial
$2$-manifolds. But, normal $d$-pseudomanifolds form a broader
class than connected combinatorial $d$-manifolds for $d \geq 3$.
In fact, any triangulation of a connected closed manifold is a
normal pseudomanifold.
%} \end{defn}

Observe that if $X$ is a normal pseudomanifold then $X$ is a
pseudomanifold. (If $\La(X)$ is not connected then, since $X$ is
connected, $\La(X)$ has two components $G_1$ and $G_2$ and two
intersecting facets $\sigma_1$, $\sigma_2$ such that $\sigma_i
\in G_i$, $i = 1, 2$. Choose $\sigma_1$, $\sigma_2$ among all
such pairs such that $\dim(\sigma_1 \cap \sigma_2)$ is maximum.
Then $\dim(\sigma_1 \cap \sigma_2) \leq d-2$ and ${\rm
lk}_X(\sigma_1 \cap \sigma_2)$ is not connected, a contradiction.)
Notice that all the links of simplices of dimensions up to $d-2$
in a normal $d$-pseudomanifold are normal pseudomanifolds.

Let $X$, $Y$ be two simplicial complexes with disjoint vertex
sets. (Since we identify isomorphic complexes, this is no real
restriction on $X$, $Y$.)  Then their {\em join} $X \ast Y$ is the
simplicial complex whose simplices are those of $X$ and of $Y$,
and the (disjoint) unions of simplices of $X$ with simplices of
$Y$. It is easy to see that if $X$ and $Y$ are combinatorial
spheres (respectively normal pseudomanifolds) then their join $X
\ast Y$ is a combinatorial sphere (respectively normal
pseudomanifold).

By a {\em subdivision} of a simplicial complex  $K$ we mean a
simplicial complex $K^{\hspace{.1mm}\prime}$ together with a
homeomorphism from $|K^{\hspace{.1mm}\prime}|$ onto $|K|$ which
is facewise linear. Two complexes $K$, $L$ have isomorphic
subdivisions if and only if $|K|$ and $|L|$ are pl homeomorphic.
Let $X$ be a pure $d$-dimensional simplicial complex and $\sigma$
be a facet of $X$, then take a symbol $v$ outside $V(X)$ and
consider the pure $d$-dimensional simplicial complex $Y$ with
vertex set $V(X) \cup \{v\}$ whose facets are facets of $X$ other
than $\sigma$ and the $(d + 1)$-sets $\tau \cup \{v\}$ where
$\tau$ runs over the $(d - 1)$-simplices in $\sigma$. Clearly,
$Y$ is a subdivision of $X$. The complex $Y$ is called the
subdivision obtained from $X$ by {\em starring a new vertex $v$ in
the facet $\sigma$}.

If $U$ is a non-empty subset of the vertex set $V(X)$ of a
simplicial complex $X$ then the simplices of $X$ which are
subsets of $U$ form a simplicial complex. This simplicial complex
is called the {\em induced subcomplex} of $X$ on the vertex set
$U$ and is denoted by $X[U]$.

\begin{defn}$\!\!\!${\bf .} \label{SC}
{\rm If $Y$ is an induced subcomplex of a simplicial complex $X$
then the {\em simplicial complement} $C(Y, X)$ of $Y$ in $X$ is
the induced subcomplex of $X$ with vertex set $V(X) \setminus
V(Y)$. By abuse of notation, for any face $\sigma$ of $X$, the
induced subcomplex of $X$ on the complement of $\sigma$ will be
denoted by $C(\sigma, X)$.}
\end{defn}

\begin{defn}$\!\!\!${\bf .} \label{DAB}
{\rm Let $\sigma_1$, $\sigma_2$ be two facets in a pure
simplicial complex $X$. Let $\psi : \sigma_1 \to \sigma_2$ be a
bijection. We shall say that $\psi$ is {\em admissible} if
($\psi$ is a bijection and) the distance between $x$ and
$\psi(x)$ in the edge graph of $X$ is $\geq 3$ for each $x\in
\sigma_1$. Notice that if $\sigma_1$, $\sigma_2$  are from
different connected components of $X$ then any bijection between
them is admissible. Also note that, in general, for the existence
of an admissible map $\psi : \sigma_1 \to \sigma_2$, the facets
$\sigma_1$ and $\sigma_2$ must be disjoint.}
\end{defn}

\begin{defn}$\!\!\!${\bf .} \label{DEHA}
{\rm Let $X$ be a weak pseudomanifold with disjoint facets
$\sigma_1$, $\sigma_2$. Let $\psi \colon \sigma_1\to \sigma_2$ be
an admissible bijection. Let $X^{\psi}$ denote the weak
pseudomanifold obtained from $X \setminus \{\sigma_1, \sigma_2\}$
by identifying $x$ with $\psi(x)$ for each $x\in \sigma_1$. Then
$X^{\psi}$ is said to be obtained from $X$ by an {\em elementary
handle addition}. If $X_1$, $X_2$ are two $d$-dimensional weak
pseudomanifolds with disjoint vertex-sets, $\sigma_i$ a facet of
$X_i$ ($i=1, 2$) and $\psi \colon \sigma_1 \to \sigma_2$ any
bijection, then $(X_1\sqcup X_2)^{\psi}$ is called an {\em
elementary connected sum} of $X_1$  and $X_2$, and is denoted by
$X_1 \#_{\psi} X_2$ (or simply by $X_1\# X_2$). Note that the
combinatorial type of $X_1 \#_{\psi} X_2$ depends on the choice
of the bijection $\psi$. However, when $X_1$, $X_2$ are connected
triangulated $d$-manifolds, $|X_1 \#_{\psi} X_2|$ is the
topological connected sum of $|X_1|$ and $|X_2|$ (taken with
appropriate orientations). Thus, $X_1 \#_{\psi} X_2$ is a
triangulated manifold whenever $X_1$, $X_2$ are triangulated
$d$-manifolds. }
\end{defn}

\begin{lemma}$\!\!\!${\bf .} \label{sd-s}
Let $N$ be a $(d- 1)$-dimensional induced subcomplex of a
$d$-dimensional simplicial complex $M$. If both $M$ and $N$ are
normal pseudomanifolds then \vspace{-1mm} \begin{enumerate}
\item[$(a)$] for any vertex $u$ of $N$ and any vertex $v$ of the
simplicial complement $C(N, M)$, there is a path $P$ $($in $M)$
joining $u$ to $v$ such that $u$ is the only vertex in $P\cap N$,
and \vspace{-1mm}
\item[$(b)$] the simplicial complement $C(N, M)$ has at most two
connected components.
\end{enumerate}
\end{lemma}

\noindent {\bf Proof.}  Part $(a)$ is trivial if $d=1$ (in which
case, $N = S^{\hspace{.2mm}0}_2$ and $M = S^{\hspace{.2mm}1}_n$).
So, assume $d > 1$ and we have the result for smaller dimensions.
Clearly, there is a path $P$ (in the edge graph of $M$) joining
$u$ to $v$ such that $P = x_1x_2 \cdots x_ky_1 \cdots y_l$ where
$x_1 = u$, $y_l = v$ and $x_i$'s are the only vertices of $P$
from $N$. Choose $k$ to be the smallest possible. We claim that $k
=1$, so that the result follows. If not, then $x_{k-1} \in {\rm
lk}_N(x_k) \subset {\rm lk}_M(x_k)$ and $y_1 \in C({\rm
lk}_N(x_k), {\rm lk}_M(x_k))$. Then, by induction hypothesis,
there is a path $Q$ in ${\rm lk}_M(x_k)$ joining $x_{k-1}$ and
$y_1$ in which $x_{k-1}$ is the only vertex from ${\rm
lk}_N(x_k)$. Replacing the part $x_{k-1}x_ky_1$ of $P$ by the
path $Q$, we get a path $P^{\hspace{.2mm} \prime}$ from $u$ to
$v$ where only the first $k-1$ vertices of $P^{\hspace{.2mm}
\prime}$ are from $N$. This contradicts the choice of $k$.

The proof of Part $(b)$ is also by induction on the dimension $d$.
The result is trivial for $d=1$. For $d > 1$, fix a vertex $u$ of
$N$. By induction hypothesis, $C({\rm lk}_N(u), {\rm lk}_M(u))$
has at most two connected components. By Part $(a)$ of this
lemma, every vertex $v$ of $C(N, M)$ is joined by a path in $C(N,
M)$ to a vertex in one of these components. Hence the result.
\hfill $\Box$

\bigskip

Let $N$ be an induced subcomplex of a simplicial complex $M$. One
says that $N$ is {\em two-sided} in $M$ if $|N|$ has a (tubular)
neighbourhood in $|M|$ homeomorphic to $|N| \times [-1, 1]$ such
that the image of $|N|$ (under this homeomorphism) is $|N|\times
\{0\}$.

\begin{lemma}$\!\!\!${\bf .} \label{m-sts}
Let $M$ be a normal pseudomanifold of dimension $d \geq 2$ and
$A$ be a set of vertices of $M$ such that the induced subcomplex
$M[A]$ of $M$ on $A$ is a $(d - 1)$-dimensional normal
pseudomanifold. Let $G$ be the graph whose vertices are the edges
of $M$ with exactly one end in $A$, two such vertices being
adjacent in $G$ if the union of the corresponding edges is a
$2$-simplex of $M$. Then $G$ has at most two connected
components. If, further, $M[A]$ is two-sided in $M$ then $G$ has
exactly two connected components.
\end{lemma}

\noindent {\bf Proof.} Let $E = V(G)$ be the set of edges of $M$
with exactly one end in $A$. For $x\in A$, set $E_x = \{e\in E :
x\in e\}$, and let $G_x = G[E_x]$ be the induced subgraph of $G$
on $E_x$. Note that $G_x$ is isomorphic to the edge graph  of
$C({\rm lk}_{M[A]}(x), {\rm lk}_M(x))$. Therefore, by Lemma
\ref{sd-s} $(b)$, $G_x$ has at most two components for each $x\in
A$. Also, for an edge $xy$ in $M[A]$, there is a $d$-simplex
$\sigma$ of $M$ such that $xy$ is in $\sigma$. Since the induced
complex $M[A]$ is $(d - 1)$-dimensional, there is a vertex $u \in
\sigma \setminus A$. Then $e_1 = xu \in E_x$ and $e_2 = yu \in
E_y$ are adjacent in $G$. Thus, if $x$, $y$ are adjacent vertices
in $M[A]$ then there is an edge of $G$ between $E_x$ and $E_y$.
Since $M[A]$ is connected and $V(G) = \cup_{x \in A} E_x$, it
follows that $G$ has at most two connected components.

Now suppose $S = M[A]$ is two-sided in $M$. Let $U$ be a tubular
neighbourhood of $|S|$ in $|M|$ such that $U \setminus |S|$ has
two components, say $U^{+}$ and $U^{-}$. Since $|S|$ is compact,
we can choose $U$ sufficiently small so that $U$ does not contain
any vertex from $V(M) \setminus A$. Then, for $e \in E$, $|e|$
meets either $U^{+}$ or $U^{-}$ but not both. Put $E^{\pm} = \{e
\in E : |e| \cap U^{\pm} \neq \emptyset\}$. Then no element of
$E^{+}$ is adjacent in $G$ with any element of $E^{-}$. From the
previous argument, one sees that each $x\in A$ is in an edge from
$E^{+}$ and in an edge from $E^{-}$. Thus, both $E^{+}$ and
$E^{-}$ are non-empty. So, $G$ is disconnected. \hfill $\Box$

\begin{lemma}$\!\!\!${\bf .} \label{LEHD}
Let $X$ be a normal $d$-pseudomanifold with an induced two-sided
standard $(d - 1)$-sphere $S$. Then there is a $d$-dimensional
weak pseudomanifold $\widetilde{X}$ such that $X$ is obtained
from $\widetilde{X}$ by elementary handle addition. Further,
\vspace{-1mm}
\begin{enumerate}
\item[$(a)$]the connected components of $\widetilde{X}$ are normal
$d$-pseudomanifolds, \vspace{-1mm}
\item[$(b)$] $\widetilde{X}$ has at most two connected components,
\vspace{-1mm}
\item[$(c)$] if $\widetilde{X}$ is not connected, then $X = Y_1 \#
Y_2$, where $Y_1$, $Y_2$ are the connected components of
$\widetilde{X}$, and \vspace{-1mm}
\item[$(d)$] if \, $C(S, X)$ is connected then $\widetilde{X}$ is
connected.
\end{enumerate}
\end{lemma}

\noindent {\bf Proof.} As above, let $E$ be the set of all edges
of $X$ with exactly one end in $S$. Let $E^{+}$ and $E^{-}$ be
the connected components of the graph $G$ (with vertex-set $E$)
defined above (cf. Lemma \ref{m-sts}). Notice that if a facet
$\sigma$ intersects $V(S)$ then $\sigma$ contains edges from $E$,
and the graph $G$ induces a connected subgraph on the set
$E_{\sigma} = \{e\in E : e\subseteq \sigma\}$. (Indeed, this
subgraph is the line graph of a complete bipartite graph.)
Consequently, either $E_{\sigma} \subseteq E^{+}$ or $E_{\sigma}
\subseteq E^{-}$. Accordingly, we say that the facet $\sigma$ is
positive or negative (relative to $S$). If a facet $\sigma$ of
$X$ does not intersect $V(S)$ then we shall say that $\sigma$ is
a neutral facet.

Let $V(S) = W$ and $V(X) \setminus V(S) = U$. Take two disjoint
sets $W^{+}$ and $W^{-}$, both disjoint from $U$, together with
two bijections $f_{\pm} \colon W \to W^{\pm}$. We define a pure
simplicial complex $\widetilde{X}$ as follows. The vertex-set of
$\widetilde{X}$ is $U \sqcup W^{+} \sqcup W^{-}$. The facets of
$\widetilde{X}$  are: (i) $W^{+}$, $W^{-}$, (ii) all the neutral
facets of $X$, (iii) for each positive facet $\sigma$ of $X$, the
set $\widetilde{\sigma} := (\sigma \cap U) \sqcup f_{+}(\sigma
\cap W)$, and (iv) for each negative facet $\tau$ of $X$, the set
$\widetilde{\tau} := (\tau \cap U) \sqcup f_{\!-}(\tau \cap W)$.
Clearly, $\widetilde{X}$ is a weak pseudomanifold. Let $\psi =
f_{-} \circ f_{+}^{-1} \colon W^{+} \to W^{-}$. It is easy to see
that $\psi$ is admissible and $X = (\widetilde{X})^{\psi}$.

Since the links of faces of dimension up to $d - 2$ in $X$ are
connected, it follows that the links of faces of dimension up to
$d - 2$ in $\widetilde{X}$ are connected. This proves $(a)$.

As $X$ is connected, choosing two vertices $f_{\pm}(x_0) \in
W^{\pm}$ of $\widetilde{X}$, one sees that each vertex of
$\widetilde{X}$ is joined by a path in the edge graph  of
$\widetilde{X}$ to either $f_{+}(x_0)$ or $f_{-}(x_0)$. Hence
$\widetilde{X}$ has at most two components. This proves $(b)$.
This arguments also shows that when $\widetilde{X}$ is
disconnected, $W^{+}$ and $W^{-}$ are facets in different
components of $\widetilde{X}$. Hence $(c)$ follows.

Observe that $C(S, X) = C(W^{+} \sqcup W^{-}, \widetilde{X})$.
Assume that $C(S, X)$ is connected. Now, for any $(d-1)$-simplex
$\tau \subseteq W^{+}$, there is a vertex $x$ in $C(S, X)$ such
that $\tau\cup\{x\}$ is a facet of $\widetilde{X}$. So, $C(S, X)$
and $W^{+}$ are in the same connected component of
$\widetilde{X}$. Similarly, $C(S, X)$ and $W^{-}$ are in the same
connected component of $\widetilde{X}$. This proves $(d)$. \hfill
$\Box$

\begin{defn}$\!\!\!${\bf .} \label{DEHD}
{\rm If $S$ is an induced two-sided $S^{\hspace{.35mm}d - 1}_{d +
1}$ in a normal $d$-pseudomanifold $X$, then the pure simplicial
complex $\widetilde{X}$ constructed above is said to be obtained
from $X$ by an} elementary handle deletion {\rm over $S$.}
\end{defn}

\begin{remark}$\!\!\!${\bf .} \label{REHD}
{\rm In Lemma \ref{LEHD}, if $X$ is a triangulated manifold then
it is easy to see that $\widetilde{X}$ is also a triangulated
manifold.}
\end{remark}

%\newpage

\section{Stacked spheres}

Let $X$ be a pure $d$-dimensional simplicial complex and $Y$ be
obtained from $X$ by starring a new vertex $v$ in a facet
$\sigma$. Clearly, $Y$ is a normal pseudomanifold if and only if
$X$ is so. Since $Y$ is a subdivision of $X$, it follows that $X$
is a combinatorial manifold (respectively, combinatorial sphere)
if and only if $Y$ is a combinatorial manifold (respectively,
combinatorial sphere). Notice that the new vertex $v$ is of
degree $d+1$ in $Y$, and when $d > 1$ the edge graph of $X$ is
the induced subgraph of the edge graph of $Y$ on the vertex set
$V(Y) \setminus \{v\}$.

Now, if $Y$ is a normal $d$-pseudomanifold, then note that for
any vertex $u$ of $Y$, ${\rm lk}_Y(u)$ is a normal $(d -
1)$-pseudomanifold, hence has at least $d + 1$ vertices. Thus,
each vertex of $Y$ has degree $\geq d + 1$. If $u$ is a vertex of
$Y$ of (minimal) degree $d+1$ and the number of vertices in $Y$
is $> d+2$, then consider the pure simplicial complex $X$ with
vertex set $V(Y) \setminus \{u\}$, whose facets are the facets of
$Y$ not passing through $u$, and the set of all $d+1$ neighbours
of $u$. We say that $X$ is obtained from $Y$ by {\em collapsing}
the vertex $u$. Clearly, this is the reverse of the operation of
starring a vertex $u$ in a facet of $X$.

\begin{defn}$\!\!\!${\bf .} \label{DSS1}
{\rm A simplicial complex $X$ is said to be a} stacked $d$-sphere
{\rm if there is a finite sequence $X_0$, $X_1, \dots, X_m$ of
simplicial complexes such that $X_0 = S^{\hspace{.2mm}d}_{d +
2}$, the standard $d$-sphere, $X_m = X$ and $X_{i}$ is obtained
from $X_{i - 1}$ by starring a new vertex in a facet of $X_{i -
1}$ for $1 \leq i \leq m$. Thus an $n$-vertex stacked $d$-sphere
is obtained from the standard $d$-sphere by $(n - d - 2)$-fold
starring. This implies that every stacked sphere is a
combinatorial sphere. Since, for $d > 1$, each starring increases
the number of edges by $d+1$, it follows that any $n$-vertex
stacked $d$-sphere has exactly ${d + 2 \choose 2} + (n - d - 2)(d
+ 1) = n(d + 1) - {d + 2 \choose 2}$ edges. }
\end{defn}

In \cite{ba2}, Barnette proved that any $n$-vertex polytopal
$d$-sphere has at least $n(d + 1) - {d + 2 \choose 2}$ edges. In
\cite{ka}, Kalai proved this result for triangulated manifolds
and also proved that, for $d \geq 3$, equality holds in this
inequality only for stacked spheres. In \cite{ta}, Tay
generalized these results to normal pseudomanifolds to prove\,:

\begin{theo}$\!\!\!${\bf .} $(${Lower Bound Theorem for
Normal Pseudomanifolds\,}$)$ \label{LBT-NPM} For $d \geq 2$, any
$n$-vertex normal $d$-pseudomanifold has at least $n(d + 1) - {d +
2 \choose 2}$ edges. For $d \geq 3$, equality holds only for
stacked spheres.
\end{theo}

In \cite{bd7}, we have presented a self-contained combinatorial
proof of Theorem \ref{LBT-NPM}. Using induction, it is not
difficult to prove the next four lemmas (see \cite{bd7} for
complete proofs).

\begin{lemma}$\!\!\!${\bf .} \label{LSS1}
Let $X$ be a normal pseudomanifold of dimension $d \geq 2$.
\vspace{-1mm}
\begin{enumerate}
\item[{$(a)$}] If $X$ is not the standard $d$-sphere then any two
vertices of degree $d + 1$ in $X$ are non-adjacent. \vspace{-1mm}
\item[{$(b)$}] If $X$ is a stacked sphere then $X$ has at least
two vertices of degree $d + 1$.
\end{enumerate}
\end{lemma}

\begin{lemma}$\!\!\!${\bf .} \label{LSS2}
Let $X$, $Y$ be normal $d$-pseudomanifolds. Suppose $Y$ is
obtained from $X$ by starring a new vertex in a facet of $X$.
Then $Y$ is a stacked sphere if and only if $X$ is a stacked
sphere.
\end{lemma}

\begin{lemma}$\!\!\!${\bf .} \label{LSS3}
The link of a vertex in a stacked sphere is a stacked sphere.
\end{lemma}

\begin{lemma}$\!\!\!${\bf .} \label{LSS4}
Any stacked sphere is uniquely determined by its edge graph.
\end{lemma}

\begin{lemma}$\!\!\!${\bf .} \label{LSS5}
Let $X_1$, $X_2$ be normal $d$-pseudomanifolds. Then $X_1 \# X_2$
is a stacked $d$-sphere if and only if both $X_1$, $X_2$ are
stacked $d$-spheres.
\end{lemma}

\noindent {\bf Proof.} Induction on the number $n \geq d + 3$ of
vertices in $X_1 \# X_2$. If $n = d + 3$ then both $X_1$, $X_2$
must be standard $d$-spheres (hence stacked spheres) and then
$X_1 \# X_2 = S^{\hspace{.2mm}0}_{2} \ast S^{\hspace{.2mm}d -
1}_{d + 1}$ is easily seen to be a stacked sphere. So, assume $n
> d + 3$, so that at least one of $X_1$, $X_2$ is not the
standard $d$-sphere. Without loss of generality, say $X_1$ is not
the standard $d$-sphere. Of course, $X = X_1 \# X_2$ is not a
standard $d$-sphere. Let $X$ be obtained from $X_1 \sqcup X_2
\setminus \{\sigma_1, \sigma_2\}$ by identifying a facet
$\sigma_1$ of $X_1$ with a facet $\sigma_2$ of $X_2$ by some
bijection. Then, $\sigma_1 = \sigma_2$ is a clique in the edge
graph of $X$, though it is not a facet of $X$. Notice that a
vertex $x \in V(X_1) \setminus \sigma_1$ is of degree $d + 1$ in
$X_1$ if and only if it is of degree $d + 1$ in $X$. If either
$X_1$ is a stacked sphere or $X$ is a stacked sphere then, by
Lemma \ref{LSS1}, such a vertex $x$ exists. Let $\widetilde{X}_1$
(respectively, $\widetilde{X}$) be obtained from $X_1$
(respectively, $X$) by collapsing this vertex $x$. Notice that
$\widetilde{X} = \widetilde{X}_1 \# X_2$. Therefore, by induction
hypothesis and Lemma \ref{LSS2}, we have: $X$ is a stacked sphere
$\Longleftrightarrow$ $\widetilde{X}$ is a stacked sphere
$\Longleftrightarrow$ both $\widetilde{X}_1$ and $X_2$ are
stacked spheres $\Longleftrightarrow$ both $X_1$ and $X_2$ are
stacked spheres.  \hfill $\Box$

\begin{defn}$\!\!\!${\bf .} \label{DSS2}
{\rm For $d \geq 2$, ${\cal K}(d)$ will denote the family of all
normal $d$-pseudomanifolds $X$ such that the link of each vertex
of $X$ is a stacked $(d - 1)$-sphere. Since all stacked spheres
are combinatorial spheres, it follows that the members of ${\cal
K}(d)$ are combinatorial $d$-manifolds. Notice that, Lemma
\ref{LSS3} says that all stacked $d$-spheres belong to the class
${\cal K}(d)$. Also, for $d \geq 2$, $K^{\hspace{.1mm}d}_{2d +
3}$  and all the simplicial complexes $K^{\hspace{.1mm}d}_{2d +
4}(p)$ constructed in Section 3 are in the class ${\cal K}(d)$
(cf. Proof of Lemma \ref{LSE2}).}
\end{defn}

\begin{lemma}$\!\!\!${\bf .} $($Walkup$)$ \label{LSS6} Let $X$ be
a normal $d$-pseudomanifold and $\psi \colon \sigma_1 \to
\sigma_2$ be an admissible bijection, where $\sigma_1, \sigma_2$
are facets of $X$. Then $X^{\psi} \in {\cal K}(d)$ if and only if
$X \in {\cal K}(d)$.
\end{lemma}

\noindent {\bf Proof.} For a vertex $v$ of $X$, let $\bar{v}$
denote the corresponding vertex of $X^{\psi}$. Observe that ${\rm
lk}_{X^{\psi}}(\bar{v})$ is isomorphic to ${\rm lk}_{X}(v)$ if $v
\in V(X) \setminus (\sigma_1 \cup \sigma_2)$ and ${\rm
lk}_{X^{\psi}}(\bar{v}) = {\rm lk}_{X}(v) \# {\rm
lk}_{X}(\psi(v))$ if $v \in \sigma_1$. The result now follows from
Lemma \ref{LSS5}. \hfill $\Box$

\begin{theo}$\!\!\!${\bf .} \label{TSS}
For $d \geq 2$, there is a unique $(3d + 4)$-vertex stacked
$d$-sphere ${\cal S} = {\cal S}^d_{3d + 4}$ which has a pair of
facets with an admissible bijection between them. Further, this
pair of facets and the admissible bijection between them is
unique up to automorphisms of ${\cal S}$.
\end{theo}

\noindent {\bf Proof.} {\sf Uniqueness\,:} Let $V^{+}$ and
$V^{-}$ be two (disjoint) facets in a $(3d + 4)$-vertex stacked
$d$-sphere ${\cal S}$, and $\psi \colon V^{+} \to V^{-}$ be an
admissible bijection. Put $V({\cal S}) = U \sqcup V^{+} \sqcup
V^{-}$. Thus, $\#(U) = d + 2$. Since $\psi$ is admissible, for
each $x \in V^{+}$, none of the $3d + 2$ vertices of ${\cal S}$
other than $x$ and $\psi(x)$ is adjacent (in the edge graph  of
${\cal S}$) with both $x$ and $\psi(x)$. Further, $x$ and
$\psi(x)$ are non-adjacent. Therefore,
\begin{equation} \label{ESS1}
\deg(x) + \deg(\psi(x)) \leq 3d + 2, ~ x \in V^{+}.
\end{equation}
Also, for $y\in U$, $y$ is adjacent to at most one vertex in the
pair $\{x, \psi(x)\}$ for each $x \in V^{+}$, and these $d + 1$
pairs partition $V({\cal S}) \setminus U$. So, each $y \in U$ has
at most $d + 1$ neighbours outside $U$. Since $y$ can have at most
$d + 1 = \#(U \setminus \{y\})$ neighbours in $U$, it follows that
 \begin{equation} \label{ESS2}
 \deg(y) \leq 2d + 2, ~ y \in U.
\end{equation}
From (\ref{ESS1}) and (\ref{ESS2}), we get by addition,
$$
\sum_{x \in V^{+}} \deg(x) + \sum_{x \in V^{+}} \deg(\psi(x)) +
\sum_{y \in U} \deg(y) \leq (d+1)(3d+2) + (d+2)(2d+2) =
(d+1)(5d+6).
$$
Now, the left hand side in this inequality is the sum of the
degrees of all the vertices of ${\cal S}$, which equals twice the
number of edges of ${\cal S}$. Thus ${\cal S}$ has at most
$(d+1)(5d+6)/2$ edges. But, as ${\cal S}$ is a $(3d+4)$-vertex
stacked $d$-sphere and $d\geq 2$, it has exactly $(3d+4)(d+1) -
{d+2 \choose 2} = (d+1)(5d+6)/2$ edges. Hence we must have
equality in (\ref{ESS1}) and (\ref{ESS2}). Thus we have equality
throughout the arguments leading to (\ref{ESS1}) and
(\ref{ESS2}). Therefore we have\,: $(a)$ $U$ is a $(d+2)$-clique
in the edge graph $G$ of ${\cal S}$, and $(b)$ for each $y\in U$
and $x\in V^{+}$, $y$ is adjacent to exactly one of the vertices
$x$ and $\psi(x)$. Notice that, since $U$, $V^{+}$ and $V^{-}$
are cliques and there is no edge between $V^{+}$ and $V^{-}$, it
follows that $G$ is completely determined by its (bipartite)
subgraph $H$ whose edges are the edges of $G$ between $U$ and
$V^{+}$.

Let $0 \leq m \leq d + 1$.

\smallskip

\noindent {\sf Claim.} There exist $x_i^{+}$, $1\leq i\leq m$, in
$V^{+}$ and $y_i$, $1\leq i\leq m$, in $U$ such that for each $i$
($1\leq i\leq m$), the $i$ vertices $y_1, \dots, y_i$ are the
only vertices from $U$ adjacent to $x_i^{+}$. Further, there is a
stacked $d$-sphere $X(m)$ with vertex-set $V({\cal S}) \setminus
\{x_i^{+} : 1\leq i\leq m\}$ whose edge graph  is the induced
subgraph $G_{m}$ of $G$ on this vertex set.

\smallskip

We prove the claim by finite induction on $m$. The claim is
trivially correct for $m = 0$ (take $X(0) = {\cal S}$, $G_0 =
G$). So, assume $1 \leq m \leq d + 1$ and the claim is valid for
all smaller values of $m$. By Lemma \ref{LSS1}, $X(m - 1)$ has at
least two vertices of degree  $d+1$ and they are non-adjacent in
$G_{m-1}$. Since each vertex of $U$ has degree $2d + 2$ in $G$, it
has degree $\geq 2d + 2 - (m - 1) > d + 1$ in $G_{m-1}$. Since
$V^{-}$ is a clique of $G_{m-1}$, at least one of the degree $d +
1$ vertices of $G_{m - 1}$ is in $V^{+} \setminus \{x^{+}_i : 1
\leq i < m\}$. Let $x^{+}_m$ be a vertex of degree $d + 1$ in
$G_{m - 1}$ from $V^{+} \setminus \{x^{+}_i : 1 \leq i < m\}$.
Notice that $x^{+}_{m - 1}$ is a vertex of degree $d + 1$ in $X(m
- 2)$; its set of neighbours in $G_{m - 2}$ is $\{y_j : 1 \leq j
\leq m - 1\} \sqcup (V^{+} \setminus \{x^{+}_i : 1 \leq i \leq m
- 1\})$. Since ${\rm lk}_{X(m - 2)}(x^{+}_{m - 1})$ is an
$S^{\hspace{.2mm}d - 1}_{d + 1}$, all the neighbours of
$x^{+}_{m-1}$ are mutually adjacent (in $G_{m-2}$ and hence) in
$G$. Thus, the vertices $y_j$, $1\leq j\leq m-1$, are adjacent in
$G$ with each vertex in $V^{+} \setminus \{x^{+}_i : 1\leq i \leq
m-1\}$. In particular, $x^{+}_m$ is adjacent (in $G$ and hence)
in $G_{m-1}$ to the $m-1$ vertices $y_j$, $1\leq j\leq m-1$ in
$U$. It is also adjacent to the $d+1-m$ vertices in $V^{+}
\setminus\{x^{+}_i : 1\leq i\leq m\}$ and to no vertex in
$V^{-}$. Since $x^{+}_m$ is of degree $d+1$ in $G_{m-1}$, it
follows that there is a unique vertex $y_m\in U\setminus \{y_i :
1\leq i\leq m-1\}$ which is adjacent to $x^{+}_m$ (in $G_{m-1}$
and hence) in $G$. By construction, $y_1, \dots, y_m$ are the
only vertices in $U$ adjacent to $x^{+}_m$. Now, let $X(m)$ be
obtained from $X(m -1)$ by collapsing the vertex $x^{+}_m$ of
degree $d+1$. By Lemma \ref{LSS2}, $X(m)$ is a stacked sphere.
Its edge graph  is the induced subgraph $G_m$ of $G$ on the
vertex-set $V({\cal S}) \setminus \{x^{+}_i : 1\leq i\leq m\}$.
This completes the induction step and hence proves the claim.

\smallskip

Now, by the final step $m = d+1$, we have named the vertices in
$V^{+}$ as $x^{+}_i$, $1\leq i\leq d+1$. We have also named $d+1$
of the vertices in $U$ as $y_i$, $1\leq i\leq d+1$. Let $y_{d+2}$
be the unique vertex in $U \setminus \{y_i : 1\leq i \leq d+1\}$.
Also, put $x_i^{-} = \psi(x^{+}_i) \in V^{-}$, $1\leq i\leq d+1$.
Thus, $x^{-}_i$ is adjacent to $y_j$ if and only if $x^{+}_i$ is
non-adjacent with $y_j$. This completes the description of the
edge graph  $G$ of ${\cal S}$. The vertices of $G$ are $x^{+}_i$,
$x^{-}_i$ ($1\leq i\leq d+1$) and $y_j$, $1\leq j\leq d+2$.
$x^{+}_i$ and $x^{+}_j$ (as well as $x^{-}_i$ and  $x^{-}_j$) are
adjacent in $G$ for $i\neq j$. $y_i$ and $y_j$ are adjacent in
$G$ for $i\neq j$. $x^{+}_i$ and $x^{-}_j$ are non-adjacent in
$G$ for all $i, j$. $x_i^{+}$ and $y_j$ are adjacent in $G$ if
and only if $j \leq i$. $x_i^{-}$ and $y_j$ are adjacent in $G$ if
and only if $j > i$.

Since the edge graph  $G$ is thus completely determined by the
given datum, Lemma \ref{LSS4} implies that ${\cal S}$ is uniquely
determined. Notice that the graph $G$ has maximum vertex degree
$2d +2$, and the set $U$ is uniquely determined by $G$ as the set
of its vertices of maximum degree. Also, the facets $V^{+}$,
$V^{-}$ are determined by $G$ as the connected components of the
induced subgraph of $G$ on the complement of $U$. Finally, the
above argument shows that the admissible bijection $\psi \colon
V^{+} \to V^{-}$ is also determined by $G$ since it must map the
unique vertex of degree $d+i$ in $V^{+}$ to the unique vertex of
degree $2d+2-i$ in $V^{-}$ ($1\leq i \leq d+1$). Notice that $S$
has an automorphism of order two which interchanges $x^{+}_i$ and
$x^{-}_{d+2-i}$ for each $i$ and interchanges $y_j$ and
$y_{d+3-j}$ for each $j$. This automorphism interchanges $V^{+}$
and $V^{-}$ and replaces $\psi$ by $\psi^{-1}$. This completes
the uniqueness proof.

\medskip

\noindent {\sf Existence of ${\cal S}^d_{3d + 4}$\,:} The
simplicial complex $\partial N^{d+1}_{3d+4}$ constructed in the
next section is a $(3d+4)$-vertex stacked $d$-sphere (cf. proof
of Lemma \ref{LSE2}) with an admissible bijection $\psi_0 \colon
B_{2d+3} \to A_{2d+3}$ (cf. the paragraph before Lemma
\ref{LSE3}). This proves the existence. \hfill $\Box$

\begin{remark}$\!\!\!${\bf .} \label{RSS}
{\rm $(a)$ The proof of Theorem \ref{TSS}, in conjunction with the
Lower Bound Theorem, actually shows the following. If $X$ is an
$n$-vertex normal $d$-pseudomanifold with an admissible
bijection, then $n \geq 3d + 4$, and equality holds only for $X =
{\cal S}^{d}_{3d + 4}$. $(b)$ If $\psi$ is the admissible
bijection on ${\cal S}^{d}_{3d + 4}$, then it is possible to
verify directly that $({\cal S}^{d}_{3d + 4})^{\psi} =
K^{\hspace{.1mm}d}_{2d + 3}$. This is also immediate from the
proof of Theorem \ref{TUK} below. }
\end{remark}

\section{Some Examples}

Recall that for any positive integer $n$, a {\em partition} of
$n$ is a finite weakly increasing sequence of positive integers
adding to $n$. The terms of the sequence are called the {\em
parts} of the partition. Let's say that a partition of $n$ is {\em
even} (respectively, {\em odd}) if it has an even (respectively,
odd) number of even parts. Let $P(n)$ (respectively $P_0(n)$,
respectively $P_1(n)$) denote the total number of partitions
(respectively even partitions, respectively odd partitions) of
$n$.

To appreciate the construction given below, it is important to
understand the growth rate of these number theoretic functions
$P_{\varepsilon}$, $\varepsilon = 0, 1$. Recall that if $f$, $g$
are two real valued functions on the set of positive integers,
then one says that $f$, $g$ are {\em asymptotically equal} (in
symbols, $f(n) \sim g(n)$) if ${\displaystyle \lim_{n\to\infty}}
\frac{f(n)}{g(n)} =1$. A famous theorem of Hardy and Ramanujan
(cf. \cite{ra}) says that
\begin{equation} \label{Pn}
P(n) \sim \frac{c_1}{n}e^{c_2\sqrt{n}} ~ \mbox{ as } ~ n\to\infty,
\end{equation}
where the absolute constants $c_1$, $c_2$ are given by
$$
c_1 = \frac{1}{4\sqrt{3}}, ~~ c_2 = \pi\sqrt{\frac{2}{3}}.
$$
We observe that\,:

\begin{lemma}$\!\!${\bf .} \label{LSE1}
$P_0(n) \sim \frac{c_1}{2n} e^{c_2 \sqrt{n}}$, $P_1(n) \sim
\frac{c_1}{2n} e^{c_2 \sqrt{n}}$ as $n \to \infty$.
\end{lemma}

\noindent {\bf Proof.} In view of (\ref{Pn}), it suffices to show
that $P_0(n) \sim \frac{1}{2} P(n)$, $P_1(n) \sim \frac{1}{2}
P(n)$ as $n \to \infty$. Now, $(p_1, \dots, p_k) \mapsto (1, p_1,
\dots, p_k)$ is a one to one function from the set of even
(respectively, odd) partitions of $n - 1$ to the set of even
(respectively, odd) partitions of $n$. Also, $(p_1, \dots, p_k)
\mapsto (p_1, \dots, p_{k - 1}, p_k + 1)$ is a one to one function
from the set of even (respectively, odd) partitions of $n-1$ to
the set of odd (respectively, even) partitions of $n$. Therefore,
$\min(P_0(n), P_1(n)) \geq \max(P_0(n - 1), P_1(n - 1))$. Since
$P_0(n - 1)+ P_1(n - 1) = P(n - 1)$, it follows that
$$
P_0(n) \geq \frac{1}{2} P(n - 1) ~ \mbox{ and } ~ P_1(n) \geq
\frac{1}{2} P(n - 1).
$$
But, from (\ref{Pn}) it follows that $P(n - 1) \sim P(n)$.
Therefore, ${\displaystyle \liminf_{n \to \infty}}
\frac{P_0(n)}{P(n)} \geq \frac{1}{2}$, ${\displaystyle \liminf_{n
\to \infty}} \frac{P_1(n)}{P(n)} \geq \frac{1}{2}$. But, $P_0(n)+
P_1(n) = P(n)$. Therefore, ${\displaystyle \lim_{n \to \infty}}
\frac{P_0(n)}{P(n)} = \frac{1}{2} = {\displaystyle \lim_{n \to
\infty}} \frac{P_1(n)}{P(n)}$. \hfill $\Box$

\bigskip

\noindent {\bf The Construction\,:} For $d \geq 2$, let $N^{\,d +
1}$ denote the pure $(d + 1)$-dimensional simplicial complex with
vertex-set $\ZZ$ (the set of all integers) such that the facets
of $N^{d + 1}$ are the sets of $d + 2$ consecutive integers. Then
$N^{d + 1}$ is a combinatorial $(d + 1)$-manifold with boundary
$M^{d} = \partial N^{d + 1}$. Now, $M^{d}$ is a combinatorial
$d$-manifold ($\in {\cal K}(d)$) and triangulates $\RR \times
S^{\hspace{.2mm}d - 1}$ (cf. \cite{ku1}). Clearly, the facets of
$M^{d}$ are of the form $\sigma_{n, i} := \{n, n + 1, \dots, n + d
+ 1\}\setminus\{n + i\}$, $1 \leq i \leq d$, $n \in \ZZ$
(intervals of length $d + 2$ minus an interior point).

For $m \geq 1$, let $N^{d + 1}_{m + d + 1}$ (respectively,
$M^{d}_{m + d +1}$) denote the induced subcomplex of $N^{\,d +
1}$ (respectively, $M^{d}$) on $m + d + 1$ consecutive vertices
(without loss of generality we may take $V(N^{d + 1}_{m + d + 1})
= V(M^{d}_{m + d + 1}) = \{1, 2, \dots, m+d+1\}$). Clearly,
$M^{d}_{m + d + 1}$ triangulates $[0, 1] \times S^{\hspace{.2mm}d
- 1}$ and $\partial M^{d}_{m + d + 1} = S^{\hspace{.2mm}d - 1}_{d
+ 1}(A_m) \sqcup S^{\hspace{.2mm}d - 1}_{d + 1}(B_m)$, where $A_m
= \{1, \dots, d + 1\}$ and $B_m = \{m + 1, \dots, m + d + 1\}$.

\begin{lemma}$\!\!\!${\bf .} \label{LSE2}
$(a)$ $\partial N^{d + 1}_{m + d + 1}$ is a stacked $d$-sphere
and $A_m$, $B_m$ are two of its facets. $(b)$ If $\psi \colon B_m
\to A_m$ is an admissible bijection then
$X^{\hspace{.1mm}d}_m(\psi) := (\partial N^{d + 1}_{m + d +
1})^{\psi}$ is a combinatorial $d$-manifold and triangulates
$S^{\hspace{.2mm}1, d - 1}(\varepsilon)$, where $\varepsilon = 0$
if $X^{\hspace{.1mm}d}_m(\psi)$ is orientable and $\varepsilon =
1$ otherwise.
\end{lemma}

\noindent {\bf Proof.} Observe that $\partial N^{d + 1}_{d + 2}$
is the standard $d$-sphere and for $m \geq 2$, $\partial N^{d +
1}_{m + d + 1}$ is obtained from $\partial N^{d + 1}_{m + d}$ by
starring the new vertex $m + d + 1$ in the facet $B_{m - 1} =
\{m, \dots, m + d\}$ of $\partial N^{d + 1}_{m + d}$. Thus,
$\partial N^{d + 1}_{m + d + 1}$ is a stacked $d$-sphere. $A_m$
is a facet of $\partial N^{d + 1}_{i + d + 1}$ for all $i\geq 1$
and from construction, $B_m$ is a facet of $\partial N^{d + 1}_{m
+ d + 1}$. This proves $(a)$.

Thus, by Lemma \ref{LSS3}, $\partial N^{d + 1}_{m + d + 1}$ is in
${\cal K}(d)$. Then, by Lemma \ref{LSS6},
$X^{\hspace{.1mm}d}_m(\psi)$ is in the class ${\cal K}(d)$. In
consequence, $X^{\hspace{.1mm}d}_m(\psi)$ is a combinatorial
$d$-manifold. Since $M^{d}_{m + d + 1}$ triangulates $[0, 1]
\times S^{\hspace{.2mm}d - 1}$ and $M^{d}_{m + d + 1} = \partial
N^{d + 1}_{m + d + 1}\setminus\{A_m, B_m\}$, it follows that
$X^{\hspace{.1mm}d}_m(\psi)$ ($= (\partial N^{d + 1}_{m + d +
1})^{\psi}$) triangulates an $S^{\hspace{.2mm}d - 1}$-bundle over
$S^{\hspace{.2mm}1}$. But, there are only two such bundles:
$S^{\hspace{.2mm}1, d - 1}(\varepsilon)$, $\varepsilon = 0, 1$
(cf. \cite[pages 134--135]{st}). This is orientable for
$\varepsilon = 0$ and non-orientable for $\varepsilon = 1$. Hence
the result. \hfill $\Box$

\bigskip

Notice that $x\in B_m$ is at a distance $\geq 3$ from $y\in A_m$
(in the edge graph of $\partial N^{d + 1}_{m + d + 1}$) if and
only if $x - y \geq 2d + 3$. Therefore, if $m \leq 2d + 2$, it is
easy to see that there is no admissible bijection $\psi \colon
B_m \to A_m$. For $m \geq 2d + 3$ the map $\psi_0 \colon B_m \to
A_m$ given by $\psi_0(m + i) = i$ is admissible. When $m = 2d +
3$, it is the only admissible map and the resulting combinatorial
manifold $X^{\,d}_{2d + 3}(\psi_0)$ is K\"{u}hnel's $K^{\,d}_{2d
+ 3}$, triangulating $S^{\hspace{.2mm}1, d - 1}(\varepsilon)$, $d
\equiv \varepsilon$ (mod 2), whose uniqueness we prove in Section
4 below. For $m \geq 2d + 3$, K\"{u}hnel and Lassmann constructed
$X^{\,d}_{m}(\psi_0)$ and proved that for $m$ odd
$X^{\,d}_{m}(\psi_0)$ is orientable if and only if $d$ is even
(cf. \cite{kl}). Here we have\,:

\begin{lemma}$\!\!\!${\bf .} \label{LSE3}
Let $m \geq 2d + 3$. If $md$ is even then for any admissible $\psi
\colon B_m \to A_m$, the combinatorial $d$-manifold
$X^{d}_{m}(\psi)$ is orientable if and only if $\psi \circ
\psi_0^{-1}$ is an even permutation. In other words, if $\psi
\circ \psi_0^{-1}$ is an even $($respectively, odd$)$ permutation
then $X^{d}_{m}(\psi)$ is a combinatorial triangulation of
$S^{\hspace{.2mm}1, d - 1}(0)$ $($respectively,
$S^{\hspace{.2mm}1, d - 1}(1))$.
\end{lemma}

\noindent {\bf Proof.} For $1 \leq k \leq m$, $1\leq i \leq d$,
let $\sigma_{k, i}$ denote the facet $\{k, k + 1, \dots, k + d +
1\} \setminus \{k + i\}$ and for $0 \leq i < j \leq d + 1$, $(i,
j) \neq (0, d+1)$, let $\sigma_{k, i, j}$ denote the $(d -
1)$-simplex $\{k, k + 1, \dots, k + d + 1\} \setminus \{k + i, k
+ j\}$ of $M^{d}_{m + d + 1}$. Consider the orientation on
$M^{d}_{m + d + 1}$ given by:
\begin{eqnarray} \label{ESE}
+ \sigma_{k, i, j} & = & (-1)^{kd + i + j} \langle k, \dots, k +
i - 1, k + i + 1, \dots, k + j - 1, k + j + 1, \dots, k + d +
1\rangle, \nonumber \\
+ \sigma_{k, i} & = & (-1)^{kd + i} \langle k, k + 1, \dots, k + i
- 1, k + i + 1, \dots, k + d + 1\rangle.
\end{eqnarray}

By an easy computation one sees that the incidence numbers
satisfy the following\,: $[\sigma_{k, i}, \sigma_{k, i, j}] = -
1$, $[\sigma_{k, j}, \sigma_{k, i, j}] = 1$ for $1\leq i < j \leq
d$, $1 \leq k \leq m$ and $[\sigma_{k, i}, \sigma_{k, 0, i}] =
1$, $[\sigma_{k + 1, i - 1}, \sigma_{k, 0, i}] = [\sigma_{k + 1,
i - 1}, \sigma_{k + 1, i - 1, d + 1}] = (- 1)^{2d -1} = -1$ for $1
\leq i \leq d$, $1 \leq k < m$. Thus, (\ref{ESE}) gives an
orientation on $M^{d}_{m + d + 1}$.

Let $\bar{\sigma}_{k, i}$ and $\bar{\sigma}_{k, i, j}$ denote the
corresponding simplices in $X^{d}_{m}(\psi_0)$. Observe that
$\bar{\sigma}_{k, 0, j} = \bar{\sigma}_{k + 1, j - 1, d + 1}$ for
$1 \leq k < m$ and $\bar{\sigma}_{m, 0, j} = \bar{\sigma}_{1,
j-1, d + 1}$. (The vertex-set of $X^{d}_{m}(\psi_0)$ is the set of
integers modulo $m$.) Then the above orientation induces an
orientation on $X^{d}_{m}(\psi_0)$. (This is well defined since
$+ \sigma_{m, 0, j} = (-1)^{md + j} \langle m + 1, \dots, m + j -
1, m + j + 1, \dots, m + d + 1 \rangle  = (- 1)^{j} \langle 1,
\dots, j - 1, j + 1, \dots, d + 1 \rangle = (- 1)^{d + (j - 1) +
(d + 1)} \langle 1, \dots, j - 1, j + 1, \dots, d + 1 \rangle = +
\sigma_{1, j - 1, d + 1}$.) Now, $[\bar{\sigma}_{m, j},
\bar{\sigma}_{m, 0, j}] = 1$, $[\bar{\sigma}_{1, j - 1},
\bar{\sigma}_{m, 0, j}] = [\bar{\sigma}_{1, j - 1},
\bar{\sigma}_{1, j - 1, d + 1}] = - 1$. Thus, $[\bar{\sigma}_{m,
j}, \bar{\sigma}_{m, 0, j}] = - [\bar{\sigma}_{1, j-1},
\bar{\sigma}_{m, 0, j}]$. Therefore, the induced orientation on
$X^{d}_{m}(\psi_0)$ is coherent. So, $X^{d}_{m}(\psi_0)$ is
orientable. This implies that $X^{d}_{m}(\psi_0)$ triangulates
$S^1\times S^{\hspace{.2mm}d - 1} = S^{\hspace{.2mm}1, d-1}(0)$.

Since $|M^{d}_{m + d + 1}|$ is homeomorphic to $|S^{\hspace{.2mm}d
- 1}_{d + 1}(B_m)| \times [0, 1]$, we can choose an orientation on
$|S^{\hspace{.2mm}d - 1}_{d + 1}(B_m)|$ so that the orientation
on $|M^{d}_{m + d + 1}|$ as the product $|S^{\hspace{.2mm}d -
1}_{d + 1}(B_m)| \times [0, 1]$ is the same as the orientation
given in (\ref{ESE}). This also induces an orientation on
$|S^{\hspace{.2mm}d - 1}_{d + 1}(A_m)|$. Let $S_B$ (respectively,
$S_A$) denote the oriented sphere $|S^{\hspace{.2mm}d - 1}_{d +
1}(B_m)|$ (respectively $|S^{\hspace{.2mm}d - 1}_{d + 1}(A_m)|$)
with this orientation. Then, as the boundary of an oriented
manifold, $\partial (|M^d_{d + m + 1}|) = S_A \cup (- S_B)$. [In
fact, it is not difficult to see that the orientation defined in
(\ref{ESE}) on $S^{\hspace{.2mm}d - 1}_{d + 1}(A_m)$
(respectively $S^{\hspace{.2mm}d - 1}_{d + 1}(B_m)$) is the same
as the orientation in $S_A$ (respectively $S_B$).]

Let $|\psi_0| \colon S_B \to S_A$ be the homeomorphism induced by
$\psi_0$. Since $|X^{d}_{m}(\psi_0)|$ is orientable, it follows
that $|\psi_0| \colon S_B \to S_A$ is orientation preserving (cf.
\cite[pages 134--135]{st}).

Therefore, $\psi \circ \psi_0^{- 1}$ is an even (respectively
odd) permutation $\Longrightarrow$ $|\psi \circ \psi_0^{- 1}|
\colon S_A \to S_A$ is orientation preserving (respectively
reversing) $\Longrightarrow$ $|\psi| = |\psi \circ \psi_0^{- 1}|
\circ |\psi_0| \colon S_B \to S_A$ is orientation preserving
(respectively reversing) $\Longrightarrow$ $|X^{d}_{m}(\psi)|$ is
orientable (respectively non-orientable). Hence, the result
follows from Lemma \ref{LSE2}. \hfill $\Box$

\bigskip

Now take $m = 2d + 4$. A bijection $\psi \colon \{2d + 5, \dots,
3d + 5\} \to \{1, \dots, d + 1\}$ is admissible for $\partial
N^{d+1}_{3d+5}$ if and only if $x - \psi(x) \geq 2d + 3$ for
$2d+5 \leq x \leq 3d+5$. It turns out that there are $2^d$
distinct admissible choices for $\psi$. But it seems difficult to
decide when two admissible choices for $\psi$ yield isomorphic
complexes $X^{\,d}_{2d + 4}(\psi)$. So, we specialize as
follows\,:

Let $p = (p_1, p_2, \dots, p_k)$ be a partition of $d + 1$. Put
$s_0 = 0$ and $s_j = \sum_{i=1}^{j} p_i$ for $1 \leq j \leq k$.
(Thus, in particular, $s_1 = p_1$ and $s_k = d + 1$.) Let $\pi_p$
be the permutation of $\{1, 2, \dots, d+1\}$ which is the product
of $k$ disjoint cycles $(s_{j - 1} + 1, s_{j - 1} +2, \dots,
s_{j})$, $1\leq j\leq k$. Notice that $\pi_p$ is an even
(respectively, odd) permutation if $p$ is an even (respectively,
odd) partition of $d + 1$. Now, define the bijection $\psi_p
\colon \{2d + 5, 2d + 6, \dots, 3d + 5\} \to \{1, 2, \dots, d +
1\}$ by $\psi_p(2d + 4 + i) = \pi_p(i)$, $1 \leq i \leq d + 1$.
Since $\pi_p(i) \leq i + 1$ for $1 \leq i \leq d + 1$, it follows
that $\psi_p$ is an admissible bijection. Clearly, the
corresponding complex $X^{d}_{2d + 4}(\psi_p)$ depends only on
the partition $p$ of $d+1$. We denote it by $K^{d}_{2d + 4}(p)$.
Note that $\pi_p = \psi_p\circ \psi_0^{-1}$. Therefore, by Lemma
\ref{LSE3}, $K^{d}_{2d + 4}(p)$ triangulates $S^{\hspace{.2mm}1,
d - 1}(0)$ (respectively, $S^{\hspace{.2mm}1, d - 1}(1)$) if $p$
is an even (respectively odd) partition of $d + 1$.

Let $G_p$ denote the non-edge graph of $K^{d}_{2d + 4}(p)$. Its
vertex-set is $V(K^{d}_{2d + 4}(p))$, and two distinct vertices
$x$, $y$ are adjacent in $G_p$ if $xy$ is not an edge of
$K^{d}_{2d + 4}(p)$. It turns out that $G_p$ has a clear
description in terms of the partition $p$. For $b\geq 1$, let
$K_{1, b}$ denote the unique graph with one vertex of degree $b$
and $b$ vertices of degree one. Also, let $p = (p_1, p_2, \dots,
p_k)$, and put $p_0 = 1$. Then a computation shows that $G_p$ is
the disjoint union of $K_{1, p_i}$, $0\leq i \leq k$. Thus, if
$p$ and $q$ are distinct partitions of $d+1$ then $G_p$ and $G_q$
are non-isomorphic (this is where our assumption that $p$, $q$
are weakly increasing sequences comes into play!) and hence
$K^{d}_{2d + 4}(p)$ and $K^{d}_{2d + 4}(q)$ are non-isomorphic
complexes. Thus we have proved\,:

\begin{theo}$\!\!\!${\bf .} \label{th1.1}
For any partition $p$ of $d + 1 \geq 3$, let $\varepsilon =
\varepsilon(p) = 0$ if $p$ is even and $= 1$ if $p$ is odd. Then
$K^{d}_{2d + 4}(p)$ is a $(2d + 4)$-vertex triangulation of
$S^{\hspace{.2mm}1, d - 1}(\varepsilon)$. Further, distinct
partitions $p$ of $d + 1$ correspond to non-isomorphic
triangulations of $S^{\hspace{.2mm}1, d - 1}(\varepsilon)$. In
consequence, for $\varepsilon = 0, 1$, there are $(2d+4)$-vertex
combinatorial triangulations of $S^{\hspace{.2mm}1, d -
1}(\varepsilon)$ and the number of non-isomorphic triangulations
is at least $P_{\varepsilon}(d+1) \sim
\frac{c_1}{2d}e^{c_2\sqrt{d}}$.
\end{theo}

This theorem provides an affirmative solution of the conjecture
(made by Lutz in \cite{lu1}) that $S^{\hspace{.2mm}1, d - 1}(1)$
can be triangulated by $2d + 4$ vertices for $d$ even. Notice that
each $(2d+4)$-vertex triangulation of $S^{\hspace{.2mm}1, d -
1}(\varepsilon)$ constructed here has $d+2$ non-edges. We
conjecture that this is the maximum possible number of non-edges.
If this is true then, for $d \equiv 1- \varepsilon$ (mod 2), our
construction yields triangulations of $S^{\hspace{.2mm}1, d -
1}(\varepsilon)$ with the minimum number of vertices and edges.

\section{Uniqueness of $K^{d}_{2d + 3}$}

Recall from Section 3 that for $d \geq 2$, $K^{d}_{2d+3}$ is the
$(2d+3)$-vertex combinatorial $d$-manifold constructed by
K\"{u}hnel in \cite{ku1}. It triangulates $S^{\hspace{.2mm}1,
d-1}(\varepsilon)$, where $\varepsilon \in \{0, 1\}$ is given by
$\varepsilon \equiv d$ (mod 2). One description of
$K^{\hspace{.1mm}d}_{2d + 3}$ is implicit in Section 3. An
equivalent (and somewhat simpler) description is as follows. It
is the boundary complex of the combinatorial $(d +1)$-manifold
with boundary whose vertices are the vertices of a cycle
$S^{\hspace{.2mm}1}_{2d + 3}$ of length $2d + 3$, and facets are
the sets of $d + 2$ vertices spanning a path in the cycle. From
this picture, it is clear that the dihedral group of order $4d +
6$ ($= {\rm Aut}(S^{\hspace{.2mm}1}_{2d + 3})$) is the full
automorphism group of $K^{d}_{2d + 3}$. Here we prove that for $d
\geq 3$, up to simplicial isomorphism, $K^{d}_{2d + 3}$ is the
unique $(2d + 3)$-vertex non-simply connected triangulated
$d$-manifold.

\begin{lemma}$\!\!\!${\bf .} $($Simplicial Alexander duality$)$
\label{LUK1} Let $L \subset L^{\prime}$ be induced subcomplexes
of a triangulated $d$-manifold $X$. Let $R\supset R^{\hspace{.2mm}
\prime}$ be the simplicial complements in $X$ of $L$ and
$L^{\prime}$ respectively. Then $H_{d - j}(L^{\prime}, L; \ZZ_2)
\cong H_{j}(R, R^{\hspace{.2mm}\prime}; \ZZ_2)$ for $0 \leq j \leq
d$.
\end{lemma}

\noindent {\bf Proof.} Fix a piecewise linear map $f \colon |X|
\to \RR$ such that for all vertices $u$ of $L$, $v$ of $R$ we have
$f(u) < f(v)$, and for all vertices $u^{\hspace{.2mm} \prime}$ of
$L^{\prime}$, $v^{\hspace{.2mm}\prime}$ of $R^{\hspace{.2mm}
\prime}$ we have $f(u^{\hspace{.2mm} \prime}) < f(v^{\hspace{.2mm}
\prime})$. Choose $c < c^{\hspace{.2mm} \prime}$ in $\RR$ such
that $f(u) < c < f(v)$ and $f(u^{\hspace{.2mm} \prime}) <
c^{\hspace{.2mm} \prime} < f(v^{\hspace{.2mm} \prime})$ for all
such $u, v, u^{\hspace{.2mm} \prime}, v^{\hspace{.2mm} \prime}$.
Define ${\cal L} = \{x \in |X| : f(x) \leq c\}$, ${\cal R} = \{x
\in |X| : f(x) > c\}$, ${\cal L}^{\prime} = \{x \in |X| : f(x)
\leq c^{\hspace{.2mm} \prime}\}$, ${\cal R}^{\hspace{.2mm}
\prime} = \{x \in |X| : f(x) > c^{\hspace{.2mm} \prime}\}$. Since
$f$ is piecewise linear, it follows that ${\cal L}, {\cal
L}^{\prime}$ are compact polyhedra (i.e., geometric carriers of
finite simplicial complexes). Also, $(|L^{\prime}|, |L|)$
(respectively $(|R|, |R^{\hspace{.2mm} \prime}|)$) is a strong
deformation retract of $({\cal L}^{\prime}, {\cal L})$
(respectively $({\cal R}, {\cal R}^{\hspace{.2mm} \prime})$).
Hence we have
\begin{eqnarray*}
&& H_{d-j}(L^{\prime}, L; \ZZ_2) \cong H_{d-j}(|L^{\prime}|, |L|;
\ZZ_2) \cong H_{d-j}({\cal L}^{\prime}, {\cal L}; \ZZ_2) \cong
H^{\hspace{.2mm}d-j}({\cal L}^{\prime}, {\cal L}; \ZZ_2) \\
&&  \cong H_{j}({\cal R}, {\cal R}^{\hspace{.2mm} \prime}; \ZZ_2)
\cong H_{j}(|R|, |R^{\hspace{.2mm} \prime}|; \ZZ_2) \cong H_{j}(R,
R^{\hspace{.2mm} \prime}; \ZZ_2) ~ \mbox { for } 0 \leq j \leq d.
\end{eqnarray*}
Here, the fourth isomorphism is because of Alexander duality (cf.
\cite[Theorem 17, Page 296]{sp}). The usual statement of this
duality refers to Alexander cohomology, but this agrees with
singular cohomology for polyhedral pairs (cf. \cite[Corollary 11,
Page 291]{sp}). Also, Alexander duality applies to orientable
closed manifolds, but any closed manifold (such as $|X|$ in our
application) is orientable over $\ZZ_2$. The third isomorphism
holds since over a field, homology and cohomology are isomorphic.
\hfill $\Box$

\begin{lemma}$\!\!\!${\bf .} \label{LUK2}
Let $X$ be a non-simply connected $n$-vertex triangulated
manifold of dimension $d \geq 3$. Then $n \geq 2d + 3$. If
further, $n = 2d + 3$, then for any facet $\sigma$ of $X$ and any
vertex $x$ outside $\sigma$, either the induced subcomplex of $X$
on $V(X) \setminus (\sigma \cup \{x\})$ is an $S^{\hspace{.2mm}d
-1}_{d + 1}$ or the induced subcomplex ${\rm lk}_X(x)[\sigma]$ of
${\rm lk}_X(x)$ on the vertex set $\sigma$ is disconnected.
\end{lemma}

\noindent {\bf Proof.} Let $\sigma$ be a facet and $C = C(\sigma,
X)$ be its simplicial complement. Choose a small (simply
connected) neighbourhood $U$ of $|\sigma|$ in $|X|$ such that $U
\cap (|X| \setminus |\sigma|)$ is homeomorphic to
$S^{\hspace{.2mm}d -1} \times (0, 1)$. Now, $|X|$ is non-simply
connected, $|X| = U \cup (|X| \setminus |\sigma|)$ and $d \geq
3$. So, by Van Kampen's theorem, $|X| \setminus |\sigma|$ is
non-simply connected. But $|C|$ is a strong deformation retract
of $|X| \setminus |\sigma|$. Therefore, $C$ is non-simply
connected.

Now fix a facet $\sigma$ of $X$. Choose an ordering $x_1, x_2,
\dots, x_n$ of $V(X)$ so that $\sigma = \{x_1, \dots, x_{d +
1}\}$. For $1 \leq i \leq n$, let $L_i$ (respectively $R_i$) be
the induced subcomplex of $X$ on the vertex-set $\{x_1, \dots,
x_i\}$ (respectively $\{x_{i + 1}, \dots, x_n\}$). Then, by Lemma
\ref{LUK1},
\begin{equation} \label{dual1}
H_{j}(R_{i}, R_{i+1}) \cong H_{d-j}(L_{i+1}, L_{i}), ~ \mbox{ for
} 0 \leq j \leq d ~ \mbox{ and } 1 \leq i < n.
\end{equation}
Here the homologies are taken with coefficients in $\ZZ_2$.

Since $L_1 = \{x_1\}$ is simply connected but $L_{n} = X$ is not,
it follows that there is a (smallest) index $i$ such that $L_i$
is simply connected but $L_{i + 1}$ is not. Note that $i \geq d +
1$. Choose this $i$. Since $L_{i + 1} = L_{i} \cup {\rm st}_{L_{i
+ 1}}(x_{i + 1})$ and $L_{i} \cap {\rm st}_{L_{i + 1}}(x_{i + 1})
= {\rm lk}_{L_{i + 1}}(x_{i + 1})$, Van Kampen's theorem implies
that ${\rm lk}_{L_{i + 1}}(x_{i + 1})$ is not connected. Hence
$H_1(L_{i + 1}, L_i) \cong H_1({\rm st}_{L_{i + 1}}(x_{i + 1}),
{\rm lk}_{L_{i + 1}}(x_{i + 1})) \cong \widetilde{H}_0({\rm
lk}_{L_{i + 1}}(x_{i + 1})) \neq \{0\}$. Thus, there is an index
$i \geq d + 1$ such that $H_1(L_{i + 1}, L_{i}) \neq \{0\}$.
Hence, from (\ref{dual1}), it follows that
\begin{equation} \label{dual2}
H_{d-2}({\rm lk}_{R_i}(x_{i+1})) \cong H_{d - 1}(R_{i},
R_{i+1})\neq \{0\} ~ \mbox{ for some } i \geq d + 1.
\end{equation}
Notice that we have $R_{i+1} \subset R_{i} \subseteq C =
C(\sigma, X)$. Since $H_{d - 1}(R_{i}, R_{i+1})\neq \{0\}$, $R_i$
contains at least two $(d-1)$-faces. Hence the number of vertices
in $R_i$ is $\geq d+1$.

First suppose $R_{i}$ has exactly $d + 1$ vertices. Since $H_{d -
2}({\rm lk}_{R_i}(x_{i + 1})) \neq \{0\}$ and ${\rm
lk}_{R_i}(x_{i + 1})$ has at most $d$ vertices, it follows that
${\rm lk}_{R_i}(x_{i + 1}) = S^{\hspace{.2mm}d - 2}_d$. Since $d
\geq 3$, it follows that $R_i$ is simply connected. As $C$ is not
simply connected, we have $R_i \subset C$ (proper inclusion).
Thus $n \geq (d + 1) + 1 + (d + 1) = 2d + 3$. Also, if the number
$n - i$ of vertices in $R_i$ is $\geq d + 2$. Then $n \geq i + d
+ 2 \geq 2d + 3$. This proves the inequality.

Now assume that $n = 2d + 3$. Let $x \not\in \sigma$ be a vertex
such that ${\rm lk}_X(x) \cap L_{d+1}$ ($= {\rm st}_X(x) \cap
L_{d+1}$) is connected. Choosing the vertex order so that  $x_{d +
2} = x$, we get that $L_{d + 2}$ is simply connected (by Van
Kampen theorem). Therefore $i \geq d + 2$. Hence $R_i$ has $\leq
n - d - 2 = d + 1$ vertices. But, $H_{d - 1}(R_{i}, R_{i + 1})
\neq \{0\}$, so that $R_{i}$ has $\geq d + 1$ vertices. Therefore
$R_i$ has exactly $d + 1$ vertices and hence $i = d + 2$. Thus,
$H_{d - 2}({\rm lk}_{R_{d + 2}}(x_{d + 3})) \cong H_{d - 1}(R_{d
+ 2}, R_{d + 3}) \neq \{0\}$. Since ${\rm lk}_{R_{d + 2}}(x_{d +
3})$ has at most $d$ vertices, it follows that ${\rm lk}_{R_{d +
2}}(x_{d + 3}) = S^{\hspace{.2mm}d - 2}_d$. Since any vertex of
$R_{d + 2}$ may be chosen to be $x_{d + 3}$ in this argument, we
get that all the vertex links of $R_{d + 2}$ are isomorphic to
$S^{\hspace{.2mm}d - 2}_d$. Hence the induced subcomplex
$R_{d+2}$ of $C$ on the vertex set $V(X) \setminus (\sigma \cup
\{x\})$ is an $S^{\hspace{.2mm}d - 1}_{d + 1}$. This proves the
lemma. \hfill $\Box$

\begin{remark}$\!\!\!${\bf .} \label{RBK}
{\rm For combinatorial manifolds, the inequality in Lemma
\ref{LUK2} is a theorem due to Brehm and K\"{u}hnel \cite{bk}.}
\end{remark}

\begin{lemma}$\!\!\!${\bf .} \label{LUK3}
Let $X$ be a $(2d + 3)$-vertex non-simply connected triangulated
manifold of dimension $d \geq 3$. Then, there is a facet $\sigma$
of $X$ such that its simplicial complement $C(\sigma, X)$
contains an induced $S^{\hspace{.2mm}d - 1}_{d + 1}$.
\end{lemma}

\noindent {\bf Proof.} Suppose the contrary. Then, by Lemma
\ref{LUK2}, for each facet $\sigma$ of $X$ and each vertex $x
\not \in \sigma$, the induced subcomplex ${\rm lk}_X(x)[\sigma]$
of ${\rm lk}_X(x)$ on $\sigma$ is disconnected. If $\tau$ were a
$(d-2)$-face of $X$ of degree 3, say with ${\rm lk}_X(\tau) =
S^{\hspace{.2mm}1}_3(\{x_1, x_2, x_3\})$, then the induced
subcomplex of ${\rm lk}_X(x_3)$ on the facet $\tau \cup \{x_1,
x_2\}$ would be connected - a contradiction. So, $X$ has no $(d -
2)$-face of degree 3. Now, no face $\gamma$ of $X$ of dimension
$e \leq d-2$ can have (minimal) degree $d - e + 1$. (In other
words, the link of $\gamma$ can not be a standard sphere.) Or
else, any $(d-2)$-face $\tau \supseteq \gamma$ of $X$ would have
degree 3. So, no standard sphere of positive dimension occurs as
a link in $X$.

Now fix a facet $\sigma$ of $X$. For each $x \in \sigma$, there
is a unique vertex $x^{\hspace{.15mm}\prime} \not \in\sigma$ such
that $(\sigma \setminus \{x\}) \cup \{x^{\hspace{.15mm}\prime}\}$
is a facet. This defines a map $x \mapsto x^{\prime}$ from
$\sigma$ to its complement. This map is injective\,: if we had
$x^{\hspace{.15mm}\prime}_1 = y = x^{\hspace{.15mm} \prime}_2$
for $x_1 \neq x_2$ then the induced subcomplex of ${\rm lk}_X(y)$
on $\sigma$ would be connected. Also, since ${\rm
lk}_X(x^{\hspace{.15mm}\prime})[\sigma]$ is disconnected, it
follows that $x$ must be an isolated vertex in ${\rm
lk}_X(x^{\hspace{.15mm}\prime})[\sigma]$. This implies that $x
x^{\hspace{.15mm}\prime}$ is an edge  of $X$, and $V({\rm lk}_X(x
x^{\hspace{.15mm}\prime})) \subseteq V(X) \setminus (\sigma
\cup\{x^{\hspace{.15mm}\prime}\})$. Hence $x x^{\hspace{.15mm}
\prime}$ is an edge of degree $\leq d + 1$. Therefore, by the
observation in the previous paragraph (with $e = 1$), $\deg_X(x
x^{\hspace{.15mm} \prime}) = d + 1$. In consequence, ${\rm lk}_X(x
x^{\hspace{.15mm} \prime})$ is a $(d + 1)$-vertex normal $(d -
2)$-pseudomanifold. But all such normal pseudomanifolds are
known\,: we must have ${\rm lk}_X(x x^{\hspace{.15mm}\prime}) =
S^{\hspace{.2mm}m}_{m + 2} \ast S^{\hspace{.2mm}n}_{n + 2}$ for
some $m, n \geq 0$ with $m + n = d-3$ (cf. \cite{bd2}). If $m >
0$ or $n > 0$ then $S^{\hspace{.2mm}1}_{3}$ occurs as a link (of
some $(d - 4)$-simplex) in this sphere and hence it occurs as the
link of a $(d - 2)$-simplex (containing $x x^{\hspace{.15mm}
\prime}$) in $X$. Hence, we must have $m = n = 0$. Thus $d = 3$
and each of the four edges $x x^{\hspace{.15mm} \prime}$ ($x\in
\sigma$) is of degree 4.

Then ${\rm lk}_X(x x^{\hspace{.15mm}\prime})$ is an
$S^{\hspace{.2mm}1}_{4} = S^{\hspace{.3mm}0}_{2} \ast
S^{\hspace{.3mm}0}_{2}$ with vertex set $V(X) \setminus
(\sigma\cup\{x^{\hspace{.15mm}\prime}\})$. In consequence,
putting $C = C(\sigma, X)$, one sees that $C$ is a 5-vertex
non-simply connected simplicial complex (by the proof of Lemma
\ref{LUK2}) such that for at least four of the vertices
$x^{\hspace{.15mm}\prime}$ in $C$, ${\rm lk}_C(x^{\hspace{.15mm}
\prime}) \supseteq S^{\hspace{.2mm}1}_{4}$. In consequence, all
${5 \choose 2} = 10$ edges occur in $C$. Since $C$ is non-simply
connected, it follows that $C$ has at least one missing triangle
(induced $S^{\hspace{.2mm}1}_{3}$), say with vertices $y_1, y_2,
y_3$. At least two of these vertices (say $y_1, y_2$) have
$S^{\hspace{.2mm}1}_{4}$ in their links. It follows that ${\rm
lk}_C(y_1) \supseteq S^{\hspace{.3mm}0}_{2}(\{y_2, y_3\}) \ast
S^{\hspace{.3mm}0}_{2}(\{y_4, y_5\})$ and ${\rm lk}_C(y_2)
\supseteq S^{\hspace{.3mm}0}_{2}(\{y_1, y_3\}) \ast
S^{\hspace{.3mm}0}_{2}(\{y_4, y_5\})$ where $y_4$, $y_5$ are the
two other vertices of $C$. Hence $C \supseteq C_0 =
(S^{\hspace{.2mm}1}_{3}(\{y_1, y_2, y_3\}) \ast
S^{\hspace{.3mm}0}_{2}(\{y_4, y_5\})) \cup \{y_4y_5\}$. But all
5-vertex simplicial complexes properly containing $C_0$ and not
containing the 2-simplex $y_1y_2y_3$ are simply connected. So, $C
= C_0$. But, then two of the vertices of $C$ (viz. $y_4$, $y_5$)
have no $S^{\hspace{.2mm}1}_{4}$ in their links, a contradiction.
This completes the proof. \hfill $\Box$

\begin{theo}$\!\!\!${\bf .} \label{TUK}
For $d \geq 3$, K\"{u}hnel's complex $K^{d}_{2d+3}$ is the only
non-simply connected $(2d+3)$-vertex triangulated manifold of
dimension $d$.
\end{theo}

\noindent {\bf Proof.} Let $X$ be a non-simply connected $(2d +
3)$-vertex triangulated manifold of dimension $d \geq 3$. By
Lemma \ref{LUK3}, $X$ must have a facet $\sigma$ such that
$C(\sigma, X)$ contains an induced subcomplex $S$ which is an
$S^{\hspace{.2mm}d - 1}_{d + 1}$. Let $x$ be the unique vertex in
$C(\sigma, X) \setminus S$. If $xy$ is a non-edge for each $y\in
\sigma$ then the normal $(d-1)$-pseudomanifold ${\rm lk}_X(x)$ is
a subcomplex of the $(d - 1)$-sphere $S$ and hence ${\rm lk}_X(x)
= S$. This implies that $C(\sigma, X)$ is the combinatorial
$d$-ball $\{x\} \ast S$. This is not possible since $C(\sigma,
X)$ is non-simply connected. Thus, $x$ forms an edge with a
vertex in $\sigma$. This implies that $C(S, X)$ is connected.

Thus, $S$ is an induced $S^{\hspace{.21mm}d - 1}_{d + 1}$ in $X$,
and $C(S, X)$ is connected. Since $d\geq 3$, $S$ is two-sided in
$X$. By Lemma \ref{LEHD}, we may delete the handle over $S$ to get
a $(3d + 4)$-vertex normal $d$-pseudomanifold $\widetilde{X}$.
Since $X$ has at most ${2d + 3 \choose 2}$ edges, it follows that
$\widetilde{X}$ has at most ${2d + 3 \choose 2} + {d + 1 \choose
2}$ edges. But ${2d + 3 \choose 2} + {d + 1 \choose 2} = (3d +
4)(d + 1) - {d + 2 \choose 2}$ is the lower bound on the number
of edges of a $(3d + 4)$-vertex normal $d$-pseudomanifold given
by the Lower Bound Theorem (cf. Theorem \ref{LBT-NPM}).
Therefore, $\widetilde{X}$ attains the lower bound, and hence, by
Theorem \ref{LBT-NPM}, $\widetilde{X}$ is a stacked sphere. Now,
Lemma \ref{LEHD} implies that $X = \widetilde{X}^{\psi}$ where
$\psi \colon \sigma_1 \to \sigma_2$ is an admissible bijection
between two facets of $\widetilde{X}$. Thus, $\widetilde{X}$ is a
$(3d + 4)$-vertex stacked $d$-sphere with an admissible bijection
$\psi$. Therefore, by Theorem \ref{TSS}, $\widetilde{X} = {\cal
S}^{\hspace{.2mm}d}_{3d + 4}$ and $\psi$ are uniquely determined,
hence so is $X = \widetilde{X}^{\psi}$. Since
$K^{\hspace{.1mm}d}_{2d + 3}$ satisfies the hypothesis, it
follows that $X = K^{\hspace{.1mm}d}_{2d + 3}$.  \hfill $\Box$

\begin{cor}$\!\!\!${\bf .} \label{CUK1}
Let $X$ be an $n$-vertex triangulation of an $S^{\hspace{.2mm}d -
1}$-bundle over $S^{\hspace{.2mm}1}$. If $d \geq 2$ then $n \geq
2d+3$. Further, if $n = 2d + 3$, then $X$ is isomorphic to
$K^{d}_{2d + 3}$.
\end{cor}

\noindent {\bf Proof.} Since an $S^{\hspace{.2mm}d - 1}$-bundle
over $S^{\hspace{.2mm}1}$ is non-simply connected, the result is
immediate from Lemma \ref{LUK2} and Theorem \ref{TUK} for $d \geq
3$. For $d = 2$, this result is classical. \hfill $\Box$

\begin{cor}$\!\!\!${\bf .} \label{CUK2}
If $d \geq 2$, $\varepsilon \equiv d$ $($mod $2)$ then
$S^{\hspace{.2mm}1, d - 1}(\varepsilon)$ has a unique $(2d +
3)$-vertex triangulation, namely $K^{d}_{2d + 3}$.
\end{cor}

\noindent {\bf Proof.} Since $S^{\hspace{.2mm}1, d -
1}(\varepsilon)$ (with $\varepsilon \equiv d$ (mod 2)) is
non-simply connected and is the geometric carrier of
$K^{\hspace{.1mm}d}_{2d + 3}$, the result is immediate from
Theorem \ref{TUK} for $d \geq 3$. For $d = 2$, this result is
classical. \hfill $\Box$

\begin{cor}$\!\!\!${\bf .} \label{CUK3}
If $d \geq 2$, $\varepsilon \not \equiv d$ $($mod $2)$ then any
triangulation of  $S^{\hspace{.2mm}1, d - 1}(\varepsilon)$
requires at least $2d + 4$ vertices. Thus, for this manifold, the
$(2d + 4)$-vertex triangulations in Section $3$ are vertex
minimal.
\end{cor}

\noindent {\bf Proof.} Since $S^{\hspace{.2mm}1, d -
1}(\varepsilon)$ (with $\varepsilon \not\equiv d$ (mod 2)) is
non-simply connected and $K^{\hspace{.1mm}d}_{2d + 3}$ does not
triangulate this space, the result is immediate from Theorem
\ref{TUK} for $d \geq 3$. For $d = 2$, this result is classical.
\hfill $\Box$

\begin{cor}$\!\!\!${\bf .} $($Walkup, Altshuler and  Steinberg$)$
\label{CUK4} $K^{3}_{9}$ is the unique $9$-vertex triangulated
$3$-manifold which is not a combinatorial $3$-sphere. In
consequence, every closed $3$-manifold other than
$S^{\hspace{.2mm}3}$ and $\TPSS$ requires at least $10$ vertices
for a triangulation.
\end{cor}

\noindent {\bf Proof.} Note that any triangulated 3-manifold is a
combinatorial 3-manifold. The result is immediate from Theorem
\ref{TUK}, since  by the Poincar\'{e}-Perelman theorem, the
$3$-sphere is the only simply connected closed $3$-manifold.
However, it is not necessary to invoke such a powerful result.
Since a simply connected 3-manifold is clearly a homology
3-sphere, and by a result of \cite{bd5} any homology 3-sphere
(other that $S^{\hspace{.2mm}3}$) requires at least 12 vertices,
the corollary follows from Theorem \ref{TUK}.  \hfill $\Box$

\bigskip

A few days after we posted the first two versions of this paper
in the arXiv (arXiv:math. GT/0610829) a similar paper \cite{css}
was posted in the arXiv (arXiv:math.CO/0611039) by Chestnut, Sapir
and Swartz. In that paper, the authors prove the uniqueness of
$K^{\hspace{.1mm}d}_{2d + 3}$ in the broader class of homology
$d$-manifolds (compared to the class of triangulated
$d$-manifolds considered here) but with a much more restrictive
topological condition (viz., $\beta_1 \neq 0$ and $\beta_2 = 0$,
compared to our hypothesis of non-simply connectedness).

\bigskip

\noindent {\bf Acknowledgement\,:} The authors thank the anonymous
referees for many useful comments which led to substantial
improvements in the presentation of this paper. The authors are
thankful to Siddhartha Gadgil for useful conversations. The
second author was partially supported by DST (Grant:
SR/S4/MS-272/05) and by UGC-SAP/DSA-IV.

\end{document}